\documentstyle[11pt,twoside]{article}
\include{graphicx}
\title{The asymptotics of a generalised Beta function}
\author{\sc R. B.\ Paris \\
{\em Division of Computing and Mathematics}, \\
{\em University of Abertay Dundee, Dundee DD1 1HG, UK}
}
\begin{document}
\def\f#1#2{\mbox{${\textstyle \frac{#1}{#2}}$}}
\def\dfrac#1#2{\displaystyle{\frac{#1}{#2}}}
\def\boldal{\mbox{\boldmath $\alpha$}}
{\newcommand{\Sgoth}{S\;\!\!\!\!\!/}
\newcommand{\bee}{\begin{equation}}
\newcommand{\ee}{\end{equation}}
\newcommand{\lam}{\lambda}
\newcommand{\ka}{\kappa}
\newcommand{\al}{\alpha}
\newcommand{\fr}{\frac{1}{2}}
\newcommand{\fs}{\f{1}{2}}
\newcommand{\g}{\Gamma}
\newcommand{\br}{\biggr}
\newcommand{\bl}{\biggl}
\newcommand{\ra}{\rightarrow}
\newcommand{\mbint}{\frac{1}{2\pi i}\int_{c-\infty i}^{c+\infty i}}
\newcommand{\mbcint}{\frac{1}{2\pi i}\int_C}
\newcommand{\mboint}{\frac{1}{2\pi i}\int_{-\infty i}^{\infty i}}
\newcommand{\gtwid}{\raisebox{-.8ex}{\mbox{$\stackrel{\textstyle >}{\sim}$}}}
\newcommand{\ltwid}{\raisebox{-.8ex}{\mbox{$\stackrel{\textstyle <}{\sim}$}}}
\renewcommand{\topfraction}{0.9}
\renewcommand{\bottomfraction}{0.9}
\renewcommand{\textfraction}{0.05}
\newcommand{\mcol}{\multicolumn}
\date{}
\maketitle
\pagestyle{myheadings}
\markboth{\hfill \sc R. B.\ Paris  \hfill}
{\hfill \sc  Asymptotics of a generalised Beta function\hfill}
\begin{abstract}
We consider the generalised Beta function introduced by Chaudhry {\it et al.\/} [J. Comp. Appl. Math. {\bf 78} (1997) 19--32] defined by
\[B(x,y;p)=\int_0^1 t^{x-1} (1-t)^{y-1} \exp \bl[\frac{-p}{4t(1-t)}\br]\,dt,\]
where $\Re (p)>0$ and the parameters $x$ and $y$ are arbitrary complex numbers. The asymptotic behaviour of $B(x,y;p)$ is obtained when (i) $p$ large, with $x$ and $y$ fixed, (ii) $x$ and $p$ large, (iii) $x$, $y$ and $p$ large and (iv) either $x$ or $y$ large, with $p$ finite. 
Numerical results are given to illustrate the accuracy of the formulas obtained.

\vspace{0.4cm}

\vspace{0.3cm}

\end{abstract}

\vspace{0.3cm}

\noindent $\,$\hrulefill $\,$

\vspace{0.2cm}

\begin{center}
{\bf 1. \  Introduction}
\end{center}
\setcounter{section}{1}
\setcounter{equation}{0}
\renewcommand{\theequation}{\arabic{section}.\arabic{equation}}
In \cite{C1}, Chaudhry {\it et al.} introduced a generalised beta function defined by the Euler-type integral\footnote{The factor 4 is introduced in the exponential for presentational convenience.}
\bee\label{e11}
B(x,y;p)=\int_0^1 t^{x-1} (1-t)^{y-1} \exp \bl[\frac{-p}{4t(1-t)}\br]\,dt,
\ee
where $\Re (p)>0$ and the parameters $x$ and $y$ are arbitrary complex numbers. When $p=0$, it is clear that when $\Re (x)>0$ and $\Re (y)>0$ the generalised function reduces to the well-known beta function $B(x,y)$ of classical analysis.
The justification for defining this extension of the beta function is given in \cite{C1} and an application of its use in defining extensions of the Gauss and confluent hypergeometric functions is discussed in \cite{C2}. It is evident from the definition in (\ref{e11}) that $B(x,y;p)$ satisfies the symmetry property
\bee\label{e12}
B(x,y;p)=B(y,x;p).
\ee

A list of useful properties of $B(x,y;p)$ is detailed by Miller in \cite{M}, where it is established that $B(x,y;p)$ may be expanded as an infinite series of Whittaker functions or Laguerre polynomials; see (\ref{a1}). He also obtained a Mellin-Barnes integral representation for $B(x,y;p)$, which we exploit in Section 2, and expressed $B(x,x\pm n;p)$ and $B(1\pm n,1;p)$, where $n$ is an integer, as finite sums of Whittaker functions. 

Our aim in this note is to derive asymptotic expansions for $B(x,y;p)$ for large $x$, $y$ and $p$. We consider (i) $|p|\ra\infty$ in $|\arg\,p|<\fs\pi$, with $x$ and $y$ fixed, (ii) $x$ and $p$ large, (iii) $x$, $y$ and $p$ large and (iv) either $x$ or $y$ large, with $p$ finite. The expansion for large $p$ is obtained using a Mellin-Barnes integral representation for $B(x,y;p)$, whereas the other cases are obtained using the method of steepest descents.

\vspace{0.6cm}

\begin{center}
{\bf 2. \ The expansion of $B(x,y;p)$ for large $p$ with $x$, $y$ finite}
\end{center}
\setcounter{section}{2}
\setcounter{equation}{0}
\renewcommand{\theequation}{\arabic{section}.\arabic{equation}}
We start with the Mellin-Barnes integral representation given by Miller \cite{M}
\bee\label{e21}
B(x,y;p)=2^{1-x-y}\pi^\fr\,\frac{1}{2\pi i}\int_{c-\infty i}^{c+\infty i} 
\frac{\g(s)\g(x+s)\g(y+s)}{\g(\fs x\!+\!\fs y\!+\!s)\g(\fs x\!+\!\fs y\!+\!\fs\!+\!s)}
p^{-s}ds
\ee
valid in $|\arg\,p|<\fs\pi$, where $c>\max \{0, -\Re (x), -\Re (y)\}$ so that the integration path lies to the right of all the poles of the integrand situated at $s=-k$, $s=-x-k$ and $s=-y-k$, $k=0, 1, 2, \ldots\,$.
Displacement of the integration path to the left over the poles followed by evaluation of the residues (assuming that no two members of the set $\{0, x, y\}$ differ by an integer -- thereby avoiding the presence of higher-order poles) yields the result that $B(x,y;p)$ can be expressed as the sum of three ${}_2F_2(-\f{1}{4}p)$ hypergeometric functions; see \cite[Eq. (1.6)]{M}.

Since there are no poles in the half-plane $\Re (s)>c$ it follows that displacement of the integration path to the right can produce no algebraic-type asymptotic expansion; see \cite[\S 5.4]{PK}. We can therefore displace the path as far to the right as we please; on such a displaced path, which we denote by $L$, the variable $|s|$ is everywhere large. The ratio of gamma functions in the integrand in (\ref{e21}) may then be expanded as an inverse factorial expansion given by \cite[p.~39, Lemma 2.2]{PK}
\[\frac{\g(s)\g(x+s)\g(y+s)}{\g(\fs x\!+\!\fs y\!+\!s)\g(\fs x\!+\!\fs y\!+\!\fs\!+\!s)}
=\sum_{j=0}^{M-1} (-)^j c_j \g(s-j-\fs)+\rho_M(s) \g(s-M-\fs),\]
where $M$ is a positive integer and $\rho_M(s)=O(1)$ as $|s|\ra\infty$ in $|\arg\,s|<\pi$. The coefficients $c_j\equiv c_j(x,y)$ are discussed below where the leading coefficient $c_0=1$. 

Substitution of the above inverse factorial expansion into the integral (\ref{e21}) then produces
\[B(x,y;p)=2^{1-x-y}\pi^\fr \bl\{\sum_{j=0}^{M-1}(-)^j c_j\ \frac{1}{2\pi i}\int_L \g(s-j-\fs)p^{-s}ds+R_M\br\},\]
where 
\[R_M=\frac{1}{2\pi i}\int_L \rho_M(s) \g(s-M-\fs)p^{-s}ds.\]
The integral may be evaluated by the well-known Cahen-Mellin integral given by (see, for example, \cite[p.~90]{PK})
\[\frac{1}{2\pi i}\int_{c-\infty i}^{c+\infty i} \g(s+\alpha)z^{-s}ds=z^\alpha e^{-z}\qquad (|\arg\,z|<\fs\pi,\ c>-\Re (\alpha))\]
to yield
\[B(x,y;p)=2^{1-x-y}\pi^\fr\bl\{p^{-\fr}e^{-p}\sum_{j=0}^{M-1} (-)^j c_j p^{-j}+R_M\br\}.\]
A bound for the remainder $R_M$ has been considered in \cite[p.~71, Lemma 2.7]{PK}, from which it follows that $R_M=O(p^{-M-\fr}e^{-p})$ as $|p|\ra\infty$ in $|\arg\,p|<\fs\pi$. 

Hence we obtain the asymptotic expansion
\bee\label{e22}
B(x,y;p)=2^{1-x-y}\pi^\fr\,p^{-\fr}e^{-p} \bl\{\sum_{j=0}^{M-1}(-)^j c_j p^{-j}+O(p^{-M})\br\}
\ee
valid as $|p|\ra\infty$ in the sector $|\arg\,p|<\fs\pi$. The expansion of $B(x,y;p)$ for large $p$ is seen to be exponentially small in $|\arg\,p|<\fs\pi$; this is a standard result when there are no poles on the right of the path in (\ref{e21}) and routine path displacement does not produce any useful asymptotic information \cite[\S 5.4]{PK}.

The coefficients $c_j$ for $j\geq 1$ can be generated by the algorithm described in \cite[\S 2.2.4]{PK}. It is found that
\[c_1=\f{1}{4}(1+x+y+2xy-x^2-y^2),\]
\[c_2=\f{1}{32}(9+6(2+xy)(x+y+xy)-(7+4xy)(x^2+y^2)-6(x^3+y^3)+x^4+y^4+14xy),\]
which are symmetrical in $x$ and $y$ as required by (\ref{e12}). A closed-form representation for $c_j$ is derived in the appendix, where it is shown that $c_j$ can be expressed in terms of a terminating ${}_3F_2(1)$ hypergeometric function given by
\bee\label{e23}
c_j\equiv c_j(x,y)=\frac{(\fs)j (y+\fs)_j}{j!}\,{}_3F_2\bl[\begin{array}{c}-j, \fs y-\fs x, \fs y-\fs x+\fs\\ \fs, y+\fs\end{array}\!; 1\br],
\ee
where $(a)_j=\g(a+j)/\g(a)$ is the Pochhammer symbol. When $x=y$, this reduces to the simpler expression
\bee\label{e24}
c_j(x,x)=\frac{(\fs)j (x+\fs)_j}{j!}.
\ee

We remark that the asymptotic expansion of $B(x,y;p)$ for $p\ra\infty$ could also have been obtained by application of the method of steepest descents, which we shall employ in the subsequent sections. See also the appendix for a different approach.
\vspace{0.6cm}

\begin{center}
{\bf 3. \ The expansion of $B(x,y;p)$ for large $x$ and $p$ with $y$ finite}
\end{center}
\setcounter{section}{3}
\setcounter{equation}{0}
\renewcommand{\theequation}{\arabic{section}.\arabic{equation}}
We consider the expansion of $B(x,y;p)$ for large $x$ and $p$, with $y$ finite,  when it is supposed that $p=ax$, where $a>0$ and $|\arg\,x|<\fs\pi$.
By the symmetry property (\ref{e12}), the same result will also cover the case of large $y$ and $p$, with $x$ finite.
From (\ref{e11}), we have 
\bee\label{e31}
B(x,y;ax)=\int_0^1 f(t) e^{-x\psi(t)}dt\qquad (|\arg\,x|<\fs\pi),
\ee
where
\[\psi(t)=\frac{a}{4t(1-t)}-\log\,t, \qquad f(t)=\frac{(1-t)^y}{t}.\]
Saddle points of the exponential factor are given by $\psi'(t)=0$; that is, at the roots of the cubic
\bee\label{e32}
t(1-t)^2+\f{1}{4}a(1-2t)=0.
\ee
We label the three saddles $t_0$, $t_1$ and $t_2$. All three saddles lie on the real axis with $t_0$ situated in the closed interval $[0,1]$, with $t_1>1$ and $t_2<0$. The $t$-plane is cut along $(-\infty,0]$. Paths of steepest descent
through the saddles $t_r$ ($r=0, 1$) are given by 
\[\Im \{e^{i\theta}(\psi(t)-\psi(t_r)\}=0,\qquad \theta=\arg\,x;\]
these paths terminate at
$t=0$ and $t=1$ in the directions $|\theta-\phi|<\fs\pi$ and $\fs\pi<\theta-\phi<\f{3}{2}\pi$, respectively, where $\phi=\arg\,t$.

When $x>0$, the integration path coincides with the steepest descent path over the saddle $t_0$; for complex $x$ in the sector $|\arg\,x|<\fs\pi$, the steepest descent path through $t_0$ becomes deformed but still terminates at $t=0$ and $t=1$; see Fig.~1. 
\begin{figure}[ht]
	\begin{center}
		\includegraphics[width=0.40\textwidth]{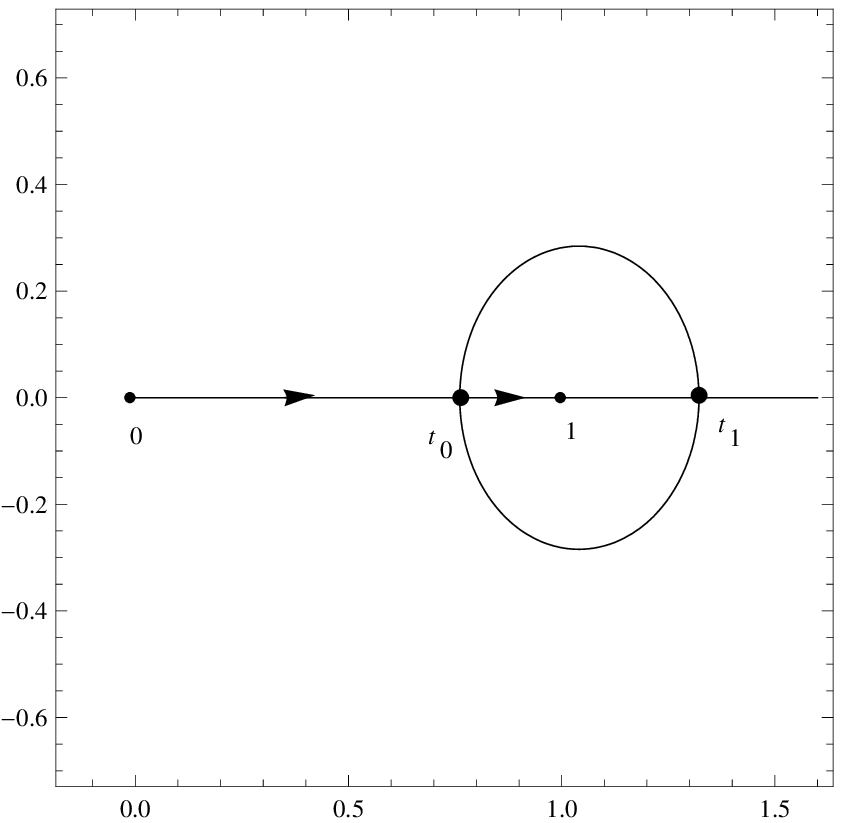}\qquad 	\includegraphics[width=0.40\textwidth]{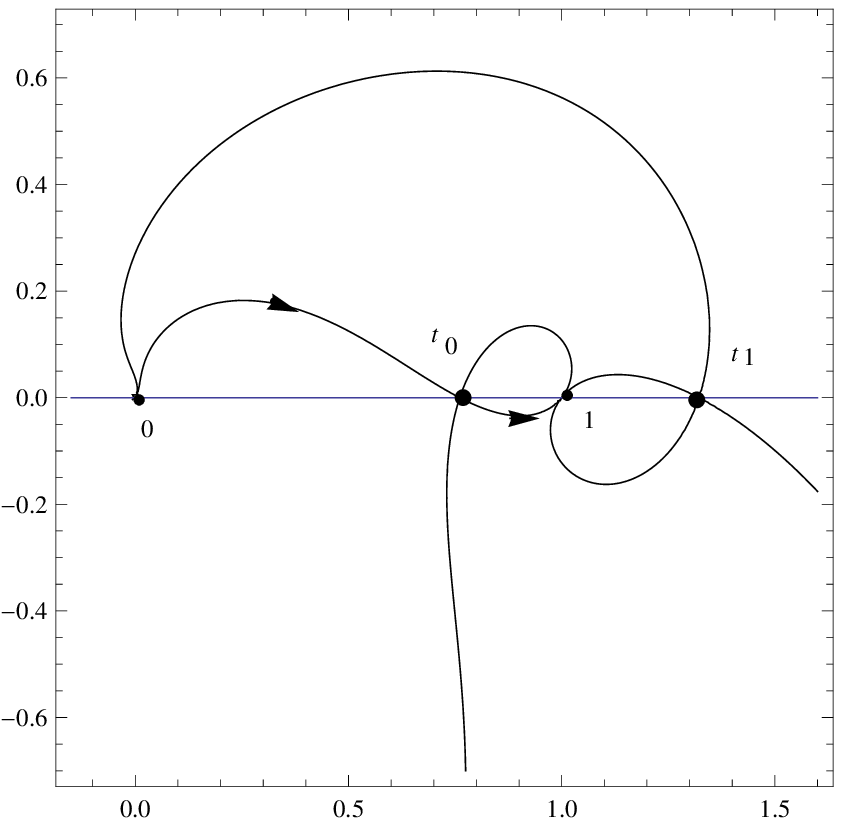}
	\caption{\small{The steepest descent and ascent paths through the saddles $t_0$ and $t_1$ (heavy dots) when $a=1/3$ and (a) $\theta=0$ and (b) $\theta=\pi/4$. The arrows indicate the integration path. In (b) the steepest ascent paths spiral round $t=0$ out to infinity passing onto adjacent Riemann surfaces. The saddle $t_2$ on the branch cut on $(-\infty,0]$ is not shown.}}\label{f1}
	\end{center}
\end{figure}
Application of the saddle-point method then yields the leading behaviour
\begin{eqnarray}
B(x,y;ax)&\sim& \sqrt{\frac{2\pi}{x\psi''(t_0)}}\,f(t_0) e^{-x\psi(t_0)}\nonumber\\
&=&\sqrt{\frac{2\pi}{x\psi''(t_0)}}\,t_0^{x-1}(1-t_0)^{y-1} \exp \bl[\frac{-ax}{4t_0(1-t_0)}\br]\label{e33}
\end{eqnarray}
as $|x|\ra\infty$ in the sector $|\arg\,x|<\fs\pi$, where some routine algebra combined with (\ref{e32}) shows that
\[\psi''(t_0)=\frac{1-3t_0+4t_0^2}{t_0^2(1-t_0)(2t_0-1)}.\]
We remark that the saddle $t_0\equiv t_0(a)$ has to be computed for a particular value of the parameter $a$, either directly from (\ref{e32}) or as a cubic root.

The asymptotic expansion of $B(x,y;ax)$ is given by \cite[p.~47]{DLMF}
\bee\label{e34}
B(x,y;ax)\sim 2e^{-x\psi(t_0)} \sum_{n=0}^\infty \frac{C_{2n} \g(n+\fs)}{x^{n+\fr}}\qquad (|x|\ra\infty,\ |\arg\,x|<\fs\pi).
\ee
The coefficients $C_{n}$ can be obtained by an inversion process and are listed for $n\leq 8$ in \cite[p.~119]{D} and for $n\leq 4$ in \cite[p.~13]{P}.
Alternatively, they can be obtained by an expansion process to yield Wojdylo's formula \cite{W} given by
\bee\label{e35}
C_{n}=\frac{1}{2a_0^{(n+1)/2}}\sum_{k=0}^{n} b_{n-k} \sum_{j=0}^k \frac{(-)^j (\fs n+\fs)_j}{j!\ a_0^j}\,B_{kj}\,;
\ee
see also \cite{LP, N}. Here $B_{kj}\equiv B_{kj}(a_1, a_2, \ldots ,a_{k-j+1})$ are the partial ordinary Bell polynomials generated by the recursion\footnote{For example, this generates the values $B_{41}=a_4$, $B_{42}=a_3^2+2a_1a_3$, $B_{43}=3a_1^2a_2$ and $B_{44}=a_1^4$.}
\[B_{kj}=\sum_{r=1}^{k-j+1} a_r B_{k-r,j-1} ,\qquad B_{k0}=\delta_{k0},\]
where $\delta_{mn}$ is the Kronecker symbol, and the coefficients $a_r$ and $b_r$ appear in the expansions
\bee\label{e36}
\psi(t)-\psi(t_0)=\sum_{r=0}^\infty a_r (t-t_0)^{r+2},\qquad f(t)=\sum_{r=0}^\infty b_r(t-t_0)^r
\ee
valid in a neighbourhood of the saddle $t=t_0$.

In numerical computations we choose a value of the parameter $a$ and compute the saddle $t_0$ from (\ref{e32}).
With a value of $y$, {\it Mathematica} is used to determine the coefficients $a_r$ and $b_r$ for $0\leq r\leq n_0$.
The coefficients $C_{2n}$ can then be calculated for $0\leq n\leq n_0$ from (\ref{e35}). We display the computed values of $C_{2n}$ for different values of $a$ and $y$ in Table 1. In Table 2, the values of the absolute relative error in the computation of $B(x,y;ax)$ from (\ref{e34}) are presented as a function of the truncation index $n$ when $x=100$.
\begin{table}[th]
\caption{\footnotesize{Values of the coefficients $C_{2n}$ (to 10dp) for different $a$ and $y$.}}
\begin{center}
\begin{tabular}{|l|l|l|l|l|}
\hline
&&&&\\[-0.3cm]
\mcol{1}{|c|}{$n$} & \mcol{1}{c|}{$a=1,\ y=1$} & \mcol{1}{c|}{$a=\fs,\ y=\f{3}{2}$} & \mcol{1}{c|}{$a=\f{3}{2},\ y=\f{5}{4}$} & \mcol{1}{c|}{$a=2,\ y=\fs$}\\
[.1cm]\hline
&&&&\\[-0.3cm]
0 & $+0.2668661228$ & $+0.1364219142$ & $+0.2036093538$ & $+0.3909054941$\\
1 & $+0.0982652355$ & $+0.2683838462$ & $+0.0762869817$ & $-0.0309094064$\\
2 & $-0.0635656655$ & $-0.1085963949$ & $-0.0456489054$ & $-0.0039290992$\\
3 & $+0.0186002666$ & $+0.0151339630$ & $+0.0137423943$ & $+0.0024209801$\\
4 & $-0.0039253710$ & $-0.0003383888$ & $-0.0026770977$ & $-0.0005115807$\\ 
5 & $+0.0012059654$ & $+0.0004533741$ & $+0.0003423270$ & $+0.0000299402$\\
[.2cm]\hline
\end{tabular}
\end{center}
\end{table}
\begin{table}[h]
\caption{\footnotesize{Values of the absolute relative error in $B(x,y;ax)$ when $x=100$ for different truncation index.}}
\begin{center}
\begin{tabular}{|l|l|l|l|l|}
\hline
&&&&\\[-0.3cm]
\mcol{1}{|c|}{$n$} & \mcol{1}{c|}{$a=1,\ y=1$} & \mcol{1}{c|}{$a=\fs,\ y=\f{3}{2}$} & \mcol{1}{c|}{$a=\f{3}{2},\ y=\f{5}{4}$} & \mcol{1}{c|}{$a=2,\ y=\fs$}\\
[.1cm]\hline
&&&&\\[-0.3cm]
0 & $1.838\times 10^{-3}$  & $9.682\times 10^{-3}$  & $1.853\times 10^{-3}$  & $3.963\times 10^{-4}$\\
1 & $1.770\times 10^{-5}$  & $5.892\times 10^{-5}$  & $1.666\times 10^{-5}$  & $7.426\times 10^{-7}$\\
2 & $1.295\times 10^{-7}$  & $2.058\times 10^{-7}$  & $1.255\times 10^{-7}$  & $1.153\times 10^{-8}$\\
3 & $9.506\times 10^{-10}$ & $1.517\times 10^{-10}$ & $8.562\times 10^{-10}$ & $8.568\times 10^{-11}$\\
4 & $1.295\times 10^{-11}$ & $9.526\times 10^{-12}$ & $5.011\times 10^{-12}$ & $2.332\times 10^{-13}$\\ 
5 & $3.688\times 10^{-12}$ & $1.933\times 10^{-13}$ & $5.472\times 10^{-14}$ & $6.917\times 10^{-15}$\\
[.2cm]\hline
\end{tabular}
\end{center}
\end{table}
\vspace{0.6cm}

\begin{center}
{\bf 4.\ The expansion of $B(x,y;p)$ for large $x$, $y$ and $p$}
\end{center}
\setcounter{section}{4}
\setcounter{equation}{0}
\renewcommand{\theequation}{\arabic{section}.\arabic{equation}}
We consider the expansion of $B(x,y;p)$ for large $x$, $y$ and $p$,  when it is supposed that $p=ax$ and $y=bx$, where $a>0$, $b>0$ and $|\arg\,x|<\fs\pi$.
From (\ref{e11}), we have
\bee\label{e41}
B(x,y;p)=\int_0^1 f(t) e^{-x\psi(t)}dt\qquad (|\arg\,x|<\fs\pi),
\ee
where
\bee\label{e41a}
\psi(t)=\frac{a}{4t(1-t)}-\log\,t-b\log (1-t),\qquad f(t)=\frac{1}{t(1-t)}.
\ee
Saddle points of the exponential factor are given by the roots of the cubic
\bee\label{e42}
t(1-t)\{1-(b+1)t\}+\f{1}{4}a(1-2t)=0.
\ee
Routine examination of this cubic shows that, when $a>0$, $b>0$, all roots are real, with one root greater than 1, one in the interval $[0,1]$ and one negative root. The distribution of the saddles is thus similar to that in Section 3, where we continue to label the saddle situated in $[0,1]$ by $t_0$.
The topology of the path of steepest descent through the saddle $t_0$, given by $\Im\{e^{i\theta}(\psi(t)-\psi(t_0)\}=0$ where $\theta=\arg\,x$, is also similar to that depicted in Fig.~1.

Accordingly, the expansion of $B(x,y;p)$ when $p=ax$ and $y=bx$, with $a>0$, $b>0$, is given by
\bee\label{e43}
B(x,bx;ax)\sim 2e^{-x\psi(t_0)} \sum_{n=0}^\infty \frac{C_{2n} \g(n+\fs)}{x^{n+\fr}}\qquad (|x|\ra\infty,\ |\arg\,x|<\fs\pi),
\ee
where the coefficients $C_{2n}$ can be determined from (\ref{e35}) when the coefficients $a_r$ and $b_r$ in (\ref{e36})
are evaluated from the definitions of $\psi(t)$ and $f(t)$ in (\ref{e41a}).

The leading behaviour is
\begin{eqnarray}
B(x,bx;ax)&\sim& \sqrt{\frac{2\pi}{x\psi''(t_0)}}\,f(t_0) e^{-x\psi(t_0)}\nonumber\\
&=&\sqrt{\frac{2\pi}{x\psi''(t_0)}}\,t_0^{x-1}(1-t_0)^{bx-1} \exp \bl[\frac{-ax}{4t_0(1-t_0)}\br]\label{e44}
\end{eqnarray}
as $|x|\ra\infty$ in the sector $|\arg\,x|<\fs\pi$, where
\[\psi''(t_0)=\frac{1-3t_0+4t_0^2}{t_0^2(1-t_0)(2t_0-1)}\bl(1-\frac{bt_0}{1-t_0}\br)+\frac{b}{t_0(1-t_0)^2}\]
and $t_0\equiv t_0(a,b)$ is the root of (\ref{e42}) situated in $t\in [0,1]$.

We note that when $b=1$ we have the result \cite{C1,M} 
\[B(x,x;p)=2^{1-2x}\pi^\fr p^{(x-1)/2} e^{-\fr p}\, W_{-\frac{1}{2}x, \frac{1}{2}x}(p)\]
in terms of the Whittaker function $W_{\kappa,\mu}(z)$; see (\ref{a1}). 
\vspace{0.6cm}

\begin{center}
{\bf 5.\ The behaviour of $B(x,y;p)$ for large $x$ and finite $y$ and $p$}
\end{center}
\setcounter{section}{5}
\setcounter{equation}{0}
\renewcommand{\theequation}{\arabic{section}.\arabic{equation}}
In this final section, we examine the behaviour of $B(x,y;p)$ for large complex $x=|x| e^{i\theta}$, with $0\leq\theta\leq\pi$, when $y$ and $p>0$ are finite. The situation when $-\pi\leq\theta\leq 0$ is analogous and, in the case of real $y$, $B(x,y;p)$ assumes conjugate values. This case has been discussed in \cite[Appendix]{C2}, but is repeated (with minor corrections) here for completeness.
By the symmetry property (\ref{e12}), the same result will also cover the case of large $y$, with $x$ and $p$ finite.

From (\ref{e11}), we have upon interchanging $x$ and $y$ (by virtue of (\ref{e12}))
\bee\label{e51}
B(x,y;p)=\int_0^1 f(t) e^{-|x|\psi(t)}dt,
\ee
where
\bee\label{e51b}
\psi(t)=\frac{\alpha}{t(1-t)}-e^{i\theta}\log (1-t), \qquad f(t)=\frac{t^{y-1}}{1-t},\qquad \alpha:=\frac{p}{4|x|}.
\ee
Because $p>0$ is a fixed parameter, the integral (\ref{e51}) is valid for arbitrary complex values of $x$ and $y$.
Saddle points of the exponential factor arise when $\psi'(t)=0$; that is, when
\bee\label{e51a}
t^2(t-1)+\alpha e^{-i\theta}(1-2t)=0.
\ee
We label the three saddles $t_0$, $t_1$ and $t_2$ as in Section 3. When $\theta=0$, all three saddles are situated on the real axis with $t_0\in [0,1]$ and $t_1>1$, $t_2<1$. As $\theta$ increases, the saddles $t_0$ and $t_2$ rotate about the origin and $t_1$ rotates about the point $t=1$. The result of this rotation is that, when $\theta=\pi$, $t_0$ and $t_2$ become a complex conjugate pair near the origin and $t_1$ is situated in the interval $[0,1]$; see Fig.~2.
\begin{figure}[ht]
	\begin{center}
		{\tiny($a$)}\includegraphics[width=0.25\textwidth]{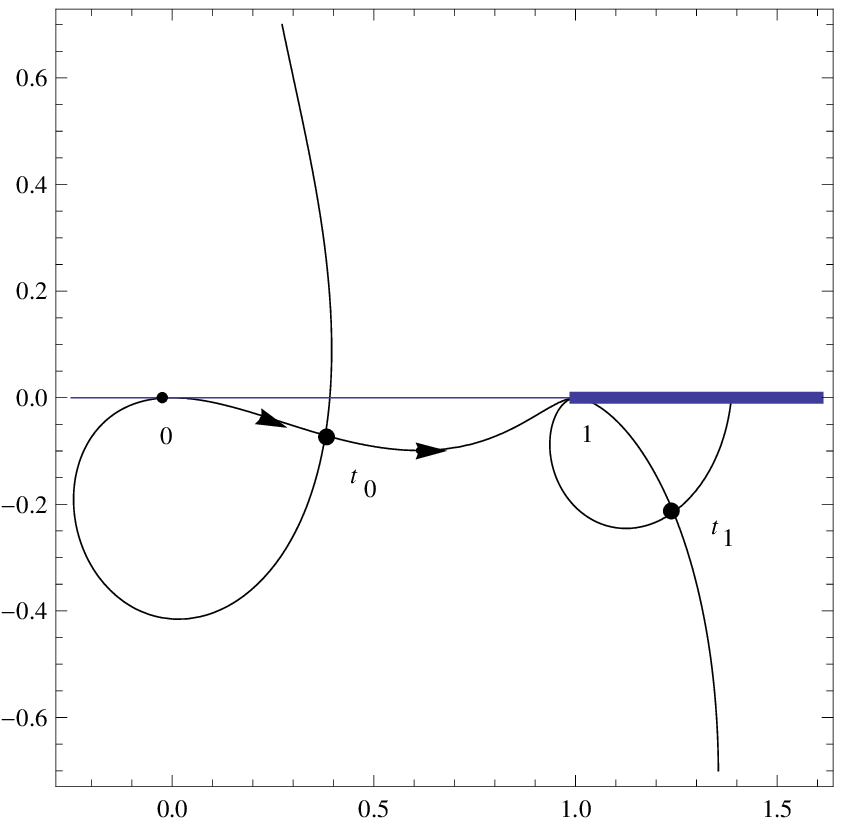}\qquad 	{\tiny($b$)}\includegraphics[width=0.25\textwidth]{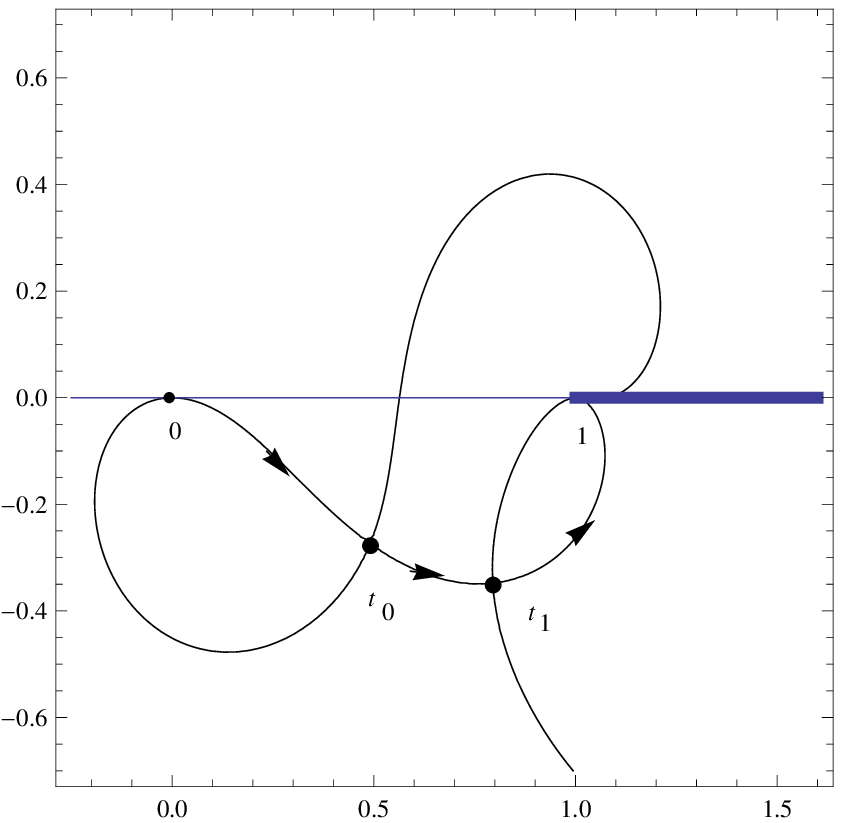}
		\qquad {\tiny($c$)}\includegraphics[width=0.25\textwidth]{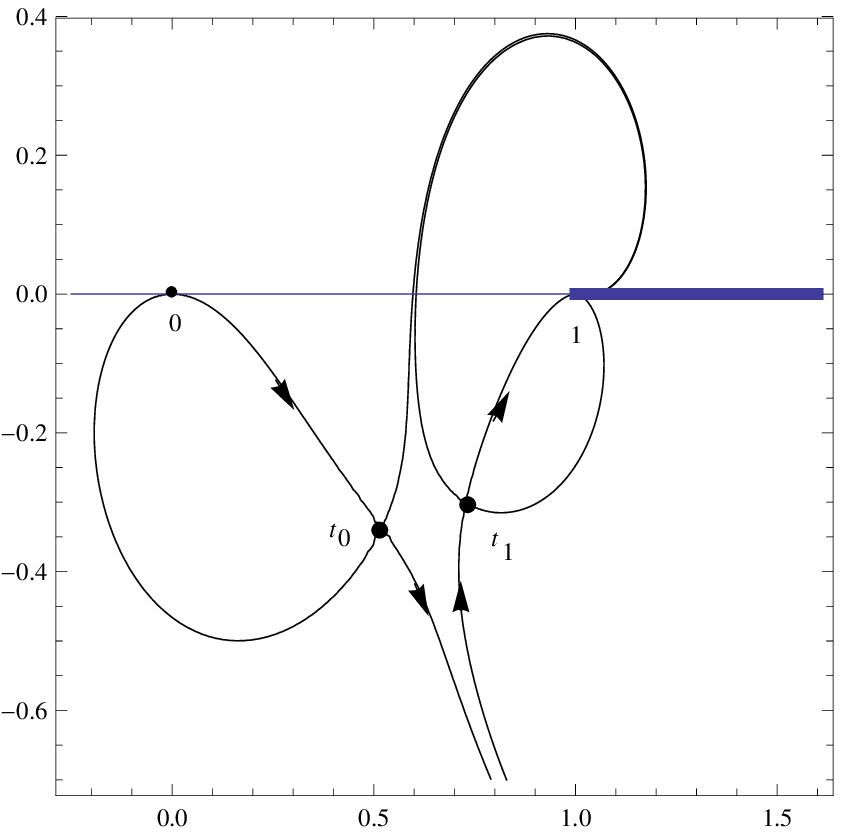}
		\vspace{0.3cm}
		
		{\tiny($d$)}\includegraphics[width=0.25\textwidth]{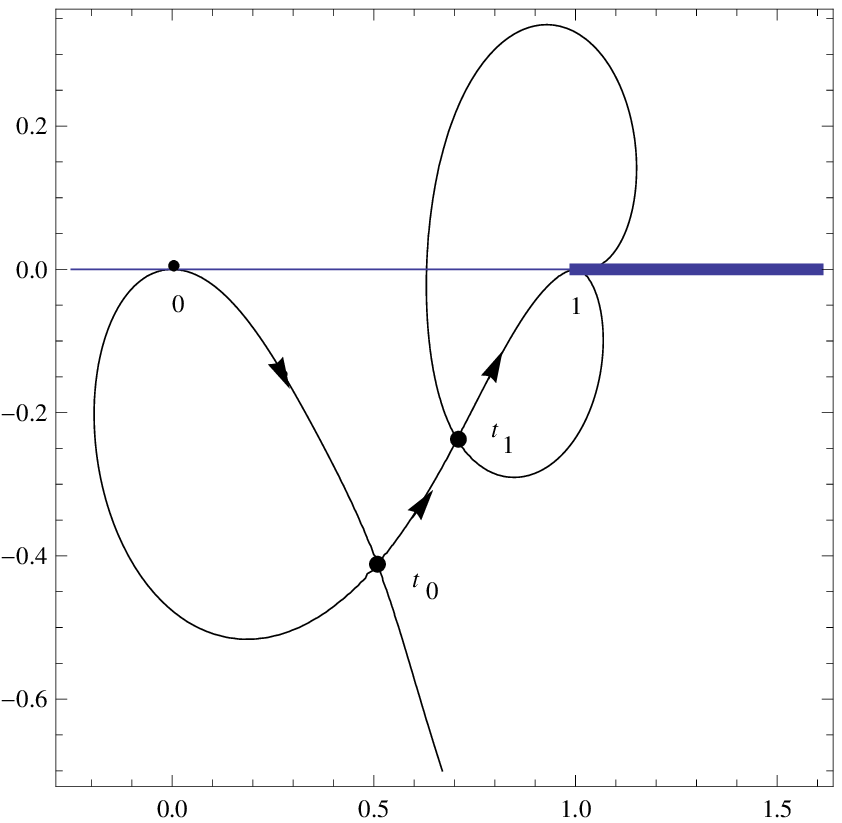}\qquad 	{\tiny($e$)}\includegraphics[width=0.25\textwidth]{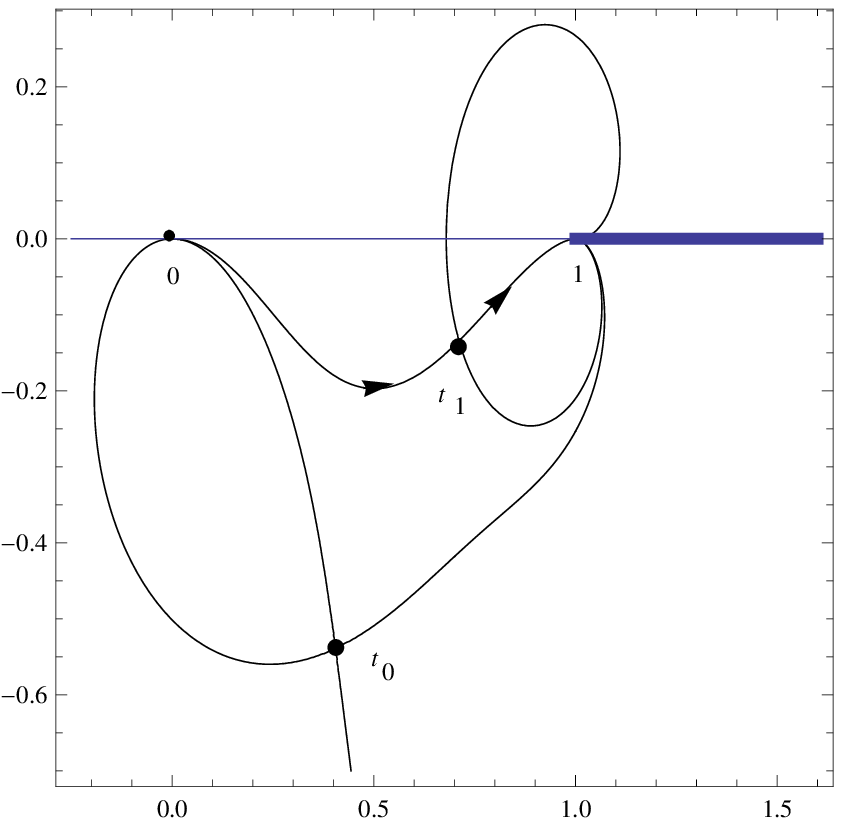}
		\qquad {\tiny($f$)}\includegraphics[width=0.25\textwidth]{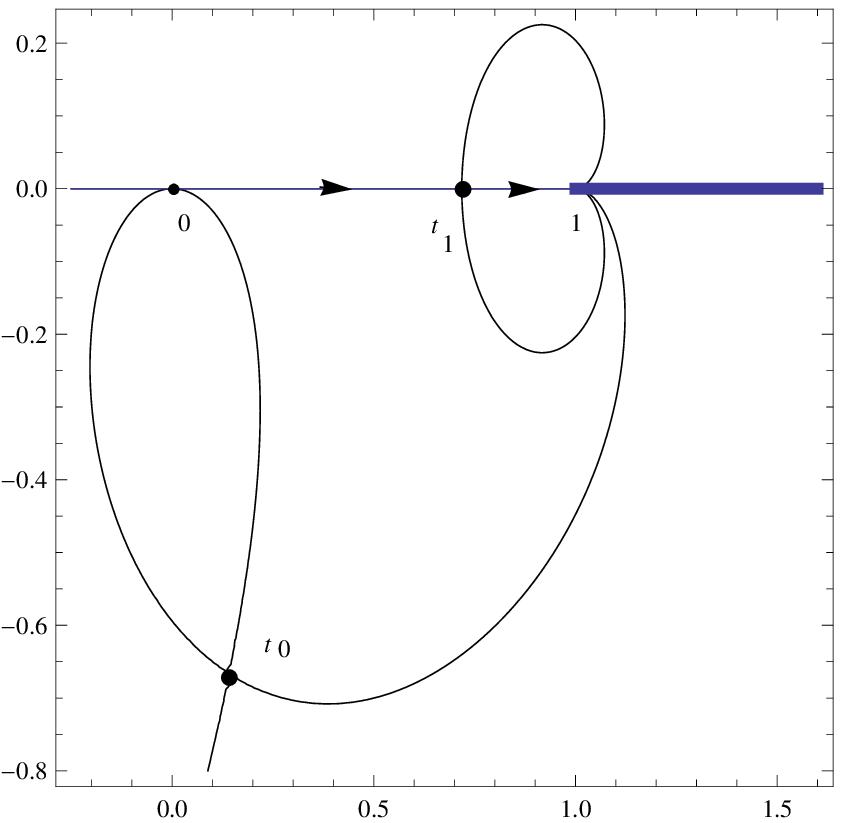}
	\caption{\small{The steepest descent and ascent paths through the saddles $t_0$ and $t_1$ (heavy dots) when $\alpha=1/3$ and (a) $\theta=0.25\pi$, (b) $\theta=\theta_0=0.65595\pi$, (c) $\theta=0.69\pi$, (d) $\theta=\theta_1=0.71782\pi$, (e) $\theta=0.80\pi$ and (f) $\theta=\pi$. The arrows indicate the integration path. The steepest ascent paths spiral round $t=1$ out to infinity passing onto adjacent Riemann surfaces. The saddle $t_2$ is not shown. The $t$-plane is cut along $[1,\infty)$.}}\label{f2}
	\end{center}
\end{figure}

When $\theta=0$, the integration path coincides with the steepest descent path passing over the saddle $t_0$ given approximately by
\[t_0\simeq\alpha^\fr-\fs\alpha\qquad (x\ra\infty).\] 
Then, with the estimates
\[x\psi(t_0)\simeq (px)^{1/2}+\f{3}{8}p,\qquad \psi''(t_0)\simeq 2\alpha^{-\fr},\]
we find by application of the saddle-point method the leading behaviour
\bee\label{e52}
B(x,y;p) \sim \sqrt{\frac{\pi}{x}}\,\bl(\frac{p}{4x}\br)^{\fr y-\frac{1}{4}} \exp \bl[-(px)^{1/2}-\frac{3}{8}p\br]\qquad (\theta=0,\ x\ra+\infty).
\ee

When $\theta=\pi$, we find from (\ref{e51a}) that the saddle $t_1$ close to the point $t=1$ is given by
\[t_1\simeq 1-\alpha+\alpha^3\qquad(|x|\ra\infty)\]
and
\[|x|\,\psi(t_1)\simeq |x|+\f{1}{4}p-|x| \log\,\alpha,\qquad \psi''(t_1)\simeq \alpha^{-2}.
\]
The integration path again coincides with the steepest descent path through $t_1$, and so we obtain the behaviour
\begin{eqnarray}
B(x,y;p)&\sim& i\sqrt{\frac{2\pi}{x}} \bl(\frac{p}{4x}\br)^{\!x}e^{\pi ix} \exp \bl[x-\frac{1}{4}p\br] \qquad (\theta=\pi,\ x\ra-\infty)\nonumber\\
&=&\sqrt{\frac{2\pi}{|x|}} \bl(\frac{p}{4|x|}\br)^{\!-|x|} \exp \bl[-|x|-\frac{1}{4}p\br].\label{e53}
\end{eqnarray}
The leading terms in (\ref{e52}) and (\ref{e53}) were given in \cite[Appendix]{C2}.

A detailed study of the topology of the steepest descent paths\footnote{The saddle $t_2$ does not enter into our consideration as it plays no role in the asymptotic evaluation of $B(x,y;p)$ when $0\leq\theta\leq\pi$.} through the saddles $t_0$ and $t_1$ when $0\leq\theta\leq\pi$ is summarised in Fig.~2 for the particular case $\alpha=\f{1}{3}$. The $t$-plane is cut along $[1,\infty)$ and paths of steepest descent either terminate at $t=0$ (with $|\arg\,t|<\fs\pi$), $t=1$ (with $|\arg (1-t)|<\fs\pi$) or at infinity. Paths that approach infinity spiral round the point $t=1$ passing onto adjacent Riemann surfaces. The figures reveal that there are two critical values of the phase $\theta$,
where the saddles $t_0$ and $t_1$ become connected (via a Stokes phenomenon). We denote these values by $\theta_0\equiv\theta_0(\alpha)$ and $\theta_1\equiv\theta_1(\alpha)$, where $\alpha$ is defined in (\ref{e51b}). The values of these critical angles are tabulated in Table 3 for different $\alpha$.

\begin{table}[h]
\caption{\footnotesize{Values of the critical angles $\theta_0$, $\theta_1$ and $\theta^*$ as a function of $\alpha=p/(4|x|)$.}}
\begin{center}
\begin{tabular}{|l|lll|}
\hline
&&&\\[-0.3cm]
\mcol{1}{|c|}{$\alpha$} & \mcol{1}{c}{$\theta_0/\pi$} & \mcol{1}{c}{$\theta_1/\pi$} & \mcol{1}{c|}{$\theta^*/\pi$}\\
[.1cm]\hline
&&&\\[-0.3cm]
0.30 & $0.603324$ & $0.752315$ & $0.688289$\\
0.25 & $0.536784$ & $0.798621$ & $0.681218$\\
0.20 & $0.476795$ & $0.840611$ & $0.672858$\\
0.15 & $0.418651$ & $0.879708$ & $0.662628$\\
0.10 & $0.358268$ & $0.916935$ & $0.649359$\\ 
0.05 & $0.288029$ & $0.953688$ & $0.629820$\\
0.01 & $0.198480$ & $0.986248$ & $0.597144$\\
[.2cm]\hline
\end{tabular}
\end{center}
\end{table}

When $0\leq\theta<\theta_0(\alpha)$, the integration path can be deformed to coincide with the steepest descent path passing over $t_0$, so that the leading behaviour in (\ref{e52}) applies in this sector. When $\theta_0(\alpha)<\theta<\theta_1(\alpha)$, the integration path is deformed to pass over both saddles $t_0$ and $t_1$, where each steepest descent path spirals out to infinity. Finally, when $\theta_1(\alpha)<\theta\leq\pi$, the integration path is deformed to pass over only the saddle $t_1$. 

Based on these considerations and on the approximation of the saddles $t_0\simeq \alpha'^\fr-\fs\alpha'$, $t_1\simeq 1+\alpha'-\alpha'^3$, where $\alpha'=p/(4x)$, the leading behaviour of $B(x,y;p)$ is found to be
\bee\label{e54}
B(x,y;p)\sim\left\{\begin{array}{ll}J_0 & 0\leq\theta<\theta_1(\alpha)\\
J_0-J_1 & \theta_1(\alpha)<\theta<\theta_2(\alpha)\\
J_1 & \theta_2(\alpha)<\theta\leq\pi \end{array}\right.
\ee
as $|x|\ra\infty$ when $0\leq\theta\leq\pi$ (with $y$ and $p>0$ finite), where
\begin{eqnarray}
J_0&:=&\sqrt{\frac{2\pi}{|x| \psi''(t_0)}} \ t_0^{y-1}(1-t_0)^{x-1} \exp \bl[\frac{-p}{4t_0(1-t_0)}\br]\label{e55}\\
&\sim&\sqrt{\frac{\pi}{x}}\,\bl(\frac{p}{4x}\br)^{\fr y-\frac{1}{4}} \exp \bl[-(px)^{1/2}-\frac{3}{8}p\br]\nonumber
\end{eqnarray}
and
\begin{eqnarray}
J_1&:=&\sqrt{\frac{2\pi}{|x| \psi''(t_1)}} \ t_1^{y-1}(1-t_1)^{x-1} \exp \bl[\frac{-p}{4t_1(1-t_1)}\br]\label{e56}\\
&\sim& i\sqrt{\frac{2\pi}{x}} \bl(\frac{p}{4x}\br)^{\!x}e^{\pi ix} \exp \bl[x-\frac{1}{4}p\br]\nonumber
\end{eqnarray}
with $\arg\,\psi''(t_r)\in [0,2\pi]$, $r=0,1$. Inspection of Table 3 shows that as $\alpha$ decreases (that is, as $|x|$ increases for fixed $p$) the angular sector $\theta_0(\alpha)\leq\theta\leq\theta_1(\alpha)$, where $B(x,y;p)$ receives a contribution from both saddles, increases. We also show in Table 3 the value of $\theta=\theta^*(\alpha)$ at which $\Re (\psi(t_0))=\Re (\psi(t_1))$ when the saddles are of the same height. We have $\theta_0(\alpha)<\theta^*(\alpha)<\theta_1(\alpha)$; then, for $\theta<\theta^*(\alpha)$ the saddle $t_0$ is dominant, whereas when $\theta>\theta^*(\alpha)$ the saddle $t_1$ is dominant in the large-$|x|$ limit.

In Table 4 we present the results of numerical calculations using the asymptotic behaviour of $B(x,y;p)$ in (\ref{e54}) compared to the values obtained by numerical integration of (\ref{e51}). The parameter values chosen correspond to $\alpha=0.01$ and the saddles $t_0$ and $t_1$ are computed from (\ref{e51a}), with the leading forms $J_0$ and $J_1$ computed from (\ref{e55}) and (\ref{e56}). It is seen from Table 3 that the exchange of dominance between the two contributory saddles arises for $\theta\simeq 0.60\pi$.
\begin{table}[h]
\caption{\footnotesize{Values of the asymptotic behaviour of $B(x,y;p)$ in (\ref{e54}) with the calculated value 
when $|x|=50$, $p=2$ ($\alpha=0.01$) and $y=\fs$ for different $\theta=\arg\,x$.}}
\begin{center}
\begin{tabular}{|l|l|l|}
\hline
&&\\[-0.3cm]
\mcol{1}{|c|}{$\theta/\pi$} & \mcol{1}{c|}{Asymptotic value} & \mcol{1}{c|}{Calculated value} \\
[.1cm]\hline
&&\\[-0.3cm]
0 & $+5.175\times 10^{-06}$ & $+5.187\times 10^{-06}$\\
0.20 & $-8.210\times 10^{-06}+2.081\times 10^{-06}i$ & $-8.223\times 10^{-06}+2.096\times 10^{-06}i$ \\
0.40 & $+3.468\times 10^{-05}-6.934\times 10^{-06}i$ & $+3.470\times 10^{-05}-7.020\times 10^{-06}i$ \\
0.50 & $+2.647\times 10^{-06}-9.853\times 10^{-05}i$ & $+2.402\times 10^{-06}-9.855\times 10^{-05}i$ \\
0.60 & $-8.837\times 10^{-04}-3.821\times 10^{-03}i$ & $-8.781\times 10^{-04}-3.823\times 10^{-03}i$ \\ 
0.70 & $-5.944\times 10^{+28}+1.659\times 10^{+28}i$ & $-5.952\times 10^{+28}+1.652\times 10^{+28}i$ \\
0.80 & $+2.786\times 10^{+54}+3.451\times 10^{+54}i$ & $+2.786\times 10^{+54}+3.459\times 10^{+54}i$ \\
1.00 & $+4.146\times 10^{+77}$  &  $+4.154\times 10^{+77}$\\
[.2cm]\hline
\end{tabular}
\end{center}
\end{table}

\vspace{0.6cm}

\begin{center}
{\bf Appendix: \ A closed-form expression for the coefficients $c_j$}
\end{center}
\setcounter{section}{1}
\setcounter{equation}{0}
\renewcommand{\theequation}{\Alph{section}.\arabic{equation}}
In this appendix we derive a closed-form expression for the coefficients $c_j$ appearing in the expansion (\ref{e22}).
Miller \cite[Eq.~(2.3a)]{M} has shown that $B(x,y;p)$ can be expressed as a convergent series of Whittaker functions
in the form
\bee\label{a1}
B(x,y;p)=2^{1-x-y}\pi^\fr\,p^{(y-1)/2}e^{-\fr p} \sum_{k=0}^\infty \frac{(\fs y\!-\!\fs x)_k(\fs\!+\!\fs y\!-\!\fs x)_k}{k!}\,W_{-k-\frac{1}{2}y, \frac{1}{2}y}(p),
\ee
where $W_{\kappa,\mu}(x)$ is the Whittaker function. For $p\ra\infty$ with bounded $k$, we have the expansion \cite[Eq.~(13.19.3)]{DLMF}
\[W_{-k-\frac{1}{2}y,\frac{1}{2}y}(p)=p^{-k-\fr}e^{-\fr p}\bl\{\sum_{n=0}^{N-1}(-)^n\frac{(\fs+k)_n(y+\fs+k)_n}{n!\,p^n}+O(p^{-N})\br\},\]
where $N$ is a positive integer. Then we obtain from (\ref{a1})
\bee\label{a2}
B(x,y;p)=2^{1-x-y}\pi^\fr\,p^{-\fr} e^{-p}\{S(x,y;p)+O(p^{-N})\},
\ee
where
\begin{eqnarray}
S(x,y;p)\!&=&\!\sum_{k=0}^{N-1}\frac{(\fs y\!-\!\fs x)_k(\fs\!+\!\fs y\!-\!\fs x)_k}{k!\,p^k}\sum_{n=0}^{N-1}(-)^n\frac{(\fs\!+\!k)_n(y\!+\!\fs\!+\!k)_n}{n!\,p^n}\nonumber\\
\!&=&\!\sum_{k=0}^{N-1}\frac{(\fs y\!-\!\fs x)_k(\fs\!+\!\fs y\!-\!\fs x)_k}{k!}\sum_{j=k}^{N-1}(-)^{j-k}\frac{(\fs\!+\!k)_{j-k}(y\!+\!\fs\!+\!k)_{j-k}}{(j-k)!\,p^j}+O(p^{-N})\nonumber
\end{eqnarray}
and we have made the change of summation index $n\ra j-k$. Use of the fact that $(-j)_k=(-)^k j!/(j-k)!$, the above double sum can be written as 
\begin{eqnarray}
&&\sum_{k=0}^{N-1}\frac{(\fs y\!-\!\fs x)_k(\fs\!+\!\fs y\!-\!\fs x)_k}{k!}\sum_{j=k}^{N-1}\frac{(-)^j}{j!}
\,\frac{\g(j\!+\!\fs)\g(y\!+\!j\!+\!\fs)}{\g(k\!+\!\fs)\g(y\!+\!\fs\!+\!k) p^j}\nonumber\\
&=&\sum_{j=0}^{N-1}\frac{(-)^j}{j!}\,\frac{(\fs)_j (y+\fs)_j}{j!\,p^j} \sum_{k=0}^j \frac{(-j)_k (\fs y\!-\!\fs x)_k(\fs\!+\!\fs y\!-\!\fs x)_k}{k!\,(\fs)_k (y+\fs)_k}\nonumber\\
&=&\sum_{j=0}^{N-1}\frac{(-)^j}{j!}\,\frac{(\fs)_j (y+\fs)_j}{j!\,p^j}\,{}_3F_2\bl[\begin{array}{c}-j, \fs y\!-\!\fs x, \fs\!+\!\fs y\!-\!\fs x\\\fs, y+\fs\end{array}\!;1\br].\label{a3}
\end{eqnarray}
upon reversal of the order of summation and identification of the inner sum over $k$ as a terminating ${}_3F_2$ series of unit argument.

Comparison of (\ref{a2}) and (\ref{a3}) with the expansion obtained in (\ref{e22}) then yields the final result
\bee\label{a4}
c_j=\frac{(\fs)j (y+\fs)_j}{j!}\,{}_3F_2\bl[\begin{array}{c}-j, \fs y\!-\!\fs x, \fs\!+\!\fs y\!-\!\fs x\\ \fs, y+\fs\end{array}\!; 1\br].
\ee

\end{document}